\documentclass[12pt]{article}

\usepackage{amssymb}
\usepackage{amsmath}
\usepackage{amsthm}
\usepackage{color}
\usepackage{verbatim}

\newtheoremstyle{plain}%
    {8pt plus2pt minus4pt}%
    {8pt plus2pt minus4pt}%
    {\itshape}%
    {}%
    {\bfseries\scshape}%
    {}%
    {6pt}
    {}%

\newtheoremstyle{remark}%
    {8pt plus2pt minus4pt}%
    {8pt plus2pt minus4pt}%
    {\upshape}
    {}%
    {\bfseries\scshape}%
    {}%
    {6pt}
    {}%

\theoremstyle{plain}
\newtheorem{theorem}{Theorem}[section]

\newtheorem{lemma}[theorem]{Lemma}
\newtheorem{proposition}[theorem]{Proposition}
\newtheorem{conjecture}[theorem]{Conjecture}
\newtheorem{problem}[theorem]{Problem}

\theoremstyle{remark}
\newtheorem{remark}[theorem]{Remark}

\newcommand{\forb}{\protect\mbox{\rm forb}}
\newcommand{\sat}{\protect\mbox{\rm sat}}
\newcommand{\Sat}{\protect\mbox{\rm SAT}}
\newcommand{\ex}{\protect\mbox{\rm ex}}
\newcommand{\di}[1]{{\dd #1}\index{#1}}
\newcommand{\C}[1]{{\protect\cal #1}}
\newcommand{\dd}{\em}
\newcommand{\claimskip}{\medskip}

\newcommand{\subs}[2]{{#1\choose #2}} 
\newcommand{\bcpf}{\noindent\emph{Proof of Claim.} }
\newcommand{\ecpf}{\qed\claimskip}

\title{On Minimum Saturated Matrices}

\author{Andrzej Dudek%
\footnote{Department of Mathematics,
Western Michigan University, Kalamazoo, MI 49008, USA}
\and
Oleg Pikhurko%
\footnote{\footnotesize {Department of Mathematical Sciences, Carnegie Mellon University, Pittsburgh, PA 15213, USA}}
\newcounter{fncmu}\setcounter{fncmu}{\value{footnote}}%
\footnote{Mathematics Institute, University of Warwick,
Coventry CV4 7AL, UK}
\thanks{\footnotesize{Partially supported by  the
National Science Foundation, Grants DMS-0758057 and DMS-1100215.}}
\and
Andrew Thomason%
\footnote{\footnotesize {Department of Pure Mathematics and Mathematical Statistics, University of Cambridge, Cambridge, CB3 0WB, UK}}
}

\begin{document}

\maketitle

\begin{abstract}

Motivated both by the work of Anstee, Griggs, and Sali on forbidden submatrices
and also by the extremal $\sat$-function for graphs, we introduce 
$\sat$-type problems for matrices. Let $\C F$ be a family of $k$-row
matrices. A matrix~$M$ is called 
\di{$\C F$-admissible} if $M$ contains no submatrix
$F\in\C F$ (as a row and column permutation of~$F$). A matrix $M$ without
repeated columns is
\di{$\C
 F$-saturated} if $M$ is $\C F$-admissible but the
addition of any column not
present in $M$ violates this property. In this paper we consider the function
$\sat(n,\C F)$ which is the \di{minimal} number of columns of an $\C
F$-saturated matrix with~$n$ rows. We establish the estimate $\sat(n,\C
F)=O(n^{k-1})$ for any family $\C F$ of $k$-row matrices and also 
compute the $\sat$-function for a few small forbidden matrices.
\end{abstract}

\section{Introduction}

First, we must introduce some simple notation. Let the shortcut `an
$n\times m$-matrix' $M$ mean a matrix with $n$ rows (which we view as
horizontal arrays) and $m$ `vertical' columns such that each entry is
0 or 1. For an $n\times m$-matrix $M$, its \di{order}
$v(M)=n$ is the number of rows and its \di{size} $e(M)=m$ is the number of
columns.  We use expressions like `an $n$-row matrix' and `an
$n$-row' to mean a matrix with $n$ rows and a row
containing $n$ elements, respectively.

For an $n\times m$-matrix $M$ and sets $A\subseteq[n]$ and
$B\subseteq[m]$, $M(A,B)$ is the $|A|\times|B|$-submatrix of $M$ formed
by the rows indexed by $A$ and the columns indexed by $B$. We use the
following obvious shorthand: $M(A,)=M(A,[m])$,
$M(A,i)=M(A,\{i\})$, etc. For example, the rows and the columns of $M$
are denoted by $M(1,),\dots,M(n,)$ and $M(,1),\dots,M(,m)$
respectively while individual entries -- by $M(i,j)$, $i\in[n]$, $j\in[m]$.

We say that a matrix $M$ is a \di{permutation} of another matrix $N$ if $M$
can be obtained from $N$ by permuting its rows and then permuting its
columns. We write $M\cong N$ in this case.  A matrix $F$ is a
\di{submatrix} of a matrix $M$ (denoted $F\subseteq M$) if we can obtain a
matrix which is a permutation of $F$ by deleting some set of rows and
columns of $M$. In other words, $F\cong M(A,B)$ for some index sets $A$ and
$B$.  The \di{transpose} of $M$ is denoted by $M^T$ (we use this notation
mostly to denote vertical columns, for typographical reasons); $(a)^i$ is
the (horizontal) sequence containing the element $a$ $i$ times. The
$n\times(m_1+m_2)$-matrix $[M_1,M_2]$ is obtained by concatenating an
$n\times m_1$-matrix $M_1$ and an $n\times m_2$-matrix $M_2$.  The
\di{complement} $1-M$ of a matrix $M$ is obtained by interchanging ones and
zeros in $M$. The \di{characteristic function} $\chi_{Y}$ of
$Y\subseteq[n]$ is the $n$-column with $i$th entry being $1$ if $i\in Y$
and $0$ otherwise.
 
Many interesting and important properties of classes of matrices can be
defined by listing forbidden submatrices. (Some authors use the term
`forbidden configurations'.) More precisely, given a family $\C F$ of
matrices (referred to as {\dd forbidden}), we say that a matrix $M$ is
\di{$\C F$-admissible} (or \di{$\C F$-free}) if $M$ contains no $F\in\C F$
as a submatrix. A \di{simple} matrix $M$ (that is, a matrix without
repeated columns) is called \di{$\C F$-saturated} (or \di{$\C F$-critical})
if $M$ is $\C F$-free but the addition of any column not present in $M$
violates this property; this is denoted by $M\in\Sat(n,\C F)$,
$n=v(M)$. Note that, although the definition requires that $M$ is simple,
we allow multiple columns in matrices belonging to $\C F$.

One well-known extremal problem is to consider $\forb(n,\C F)$, the
maximal size of a simple $\C F$-free matrix with $n$ rows or,
equivalently, the maximal size of $M\in\Sat(n,\C F)$.  Many different
results on the topic have been obtained; we refer the reader to a recent survey
by Anstee~\cite{anstee:survey}. We just want to mention a remarkable fact
that one of the first $\forb$-type results, namely formula~\eqref{eq::sauer}
here, proved independently by Vapnik and
Chervonenkis~\cite{vapnik+chervonenkis:71}, Perles and Shelah~\cite{shelah:72}, and Sauer~\cite{sauer:73}, was motivated by
such different topics as probability, logic, and a problem of Erd\H os
on infinite set systems. 

The $\forb$-problem is reminiscent of the Tur\'an function $\ex(n,\C F)$:
given a family $\C F$ of forbidden graphs, $\ex(n,\C F)$ is the maximal
size of an $\C F$-free graph on $n$ vertices not containing any member of
$\C F$ as a subgraph (see \di{e.g.}\
surveys~\cite{furedi:91,sidorenko:95,keevash:11}). Erd\H{o}s, Hajnal, and
Moon~\cite{erdos+hajnal+moon:64} considered the `dual' function $\sat(n,\C
F)$, the \emph{minimal} size of a maximal $\C F$-free graph on $n$
vertices.  This is an active area of extremal graph theory; see the dynamic
survey by Faudree, Faudree, and Schmitt~\cite{faudree+faudree+schmitt:DS}.

Here we consider the `dual' of the
$\forb$-problem for matrices. Namely, we are interested in the value of
$\sat(n,\C F)$, the {\em minimal} size of an $\C F$-saturated matrix
with $n$ rows:$$
 \sat(n,\C F)=\min\{e(M) : M\in \Sat(n,\C F)\}.
 $$
 We decided to use the same notation as for its graph counterpart. This
 should not cause any confusion as this paper will deal with matrices.
 Obviously, $\sat(n,\C F)\le\forb(n,\C F)$. If $\C F=\{F\}$ consists of a
 single forbidden matrix $F$ then we write $\Sat(n,F)=\Sat(n,\{F\})$, and
 so on.

We denote by $T_k^l$ the simple $k\times\binom kl$-matrix consisting
of all $k$-columns with exactly $l$ ones and by $K_k$ -- the $k\times
2^k$ matrix of all possible columns of order $k$. Naturally, $T_k^{\le
l}$ denotes the $k\times f(k,l)$-matrix consisting of all distinct
columns with at most $l$ ones, and so on, where we use the shortcut 
 $$
 f(k,l)= {k\choose 0}+{k\choose 1}+\dots+{k\choose l}.
 $$
Vapnik and Chervonenkis~\cite{vapnik+chervonenkis:71},
Perles and Shelah~\cite{shelah:72}, and
Sauer~\cite{sauer:73} showed independently that
\begin{equation}\label{eq::sauer}
 \forb(n,K_k)=f(n,k-1).
\end{equation}
Formula~\eqref{eq::sauer} turns out to play a significant role in our study.

This paper is organized as follows. In~\S\ref{hsfc} we give some general
results about the $\sat$-function, the principal one being
Theorem~\ref{th::hsf:k-1} which states that $\sat(n,\C F)=O(n^{k-1})$ holds
for any family $\C F$ of $k$-row matrices. Turning to specific matrices,
in~\S\ref{hsfcases} we compute $\sat(n,K_k)$ for $k=2$ and $k=3$. By
Theorem~\ref{th::hsf:k-1}, $\sat(n,K_2)$ can grow at most linearly, and
indeed it is linear in~$n$. Surprisingly, though, $\sat(n,K_3)$ is constant
for $n\ge4$. Finally, in~\S\ref{smalls}, we examine a selection of small
matrices $F$ to see how $\sat(n,F)$ behaves. In particular, we find some
$F$ for which the function grows and other $F$ for which it is constant (or
bounded): it would be interesting to determine a criterion for when
$\sat(n,F)$ is bounded, but we cannot guess one from the present data.

\section{General Results}\label{hsfc}

Here we present some results dealing with $\sat(n,\C F)$
for a general family $\C F$.

The following simple observation can be useful in tackling these
problems. Let $M'$ be obtained from $M\in\Sat(n,\C F)$ by {\dd
  duplicating}\index{row duplication} the $n$th row of $M$, that is, we let
$M'([n],)=M$ and $M'(n+1,)=M(n,)$. Suppose that $M'$ is $\C
F$-admissible. Complete $M'$, by adding columns in an arbitrary way, to an
$\C F$-saturated matrix.  Let $C$ be any added $(n+1)$-column. As both
$M'([n],)$ and $M'([n-1]\cup\{n+1\},)$ are equal to $M\in\Sat(n,\C F)$, we
conclude that both $C([n])$ and $C([n-1]\cup\{n+1\})$ must be columns of
$M$. As $C$ is not an $M'$-column, $C=(C',b,1-b)$ where $b\in\{0,1\}$ and
$C'$ is some $(n-1)$-column such that both $(C',0)$ and $(C',1)$ are
columns of $M$. This implies that $\sat(n+1,\C F)\le e(M)+2d$, where $d$ is
the number of pairs of equal columns in $M$ after we delete the $n$th
row. In particular, the following theorem follows.

\begin{theorem}\label{th::hsf:const} Suppose that $F$ is a matrix with no
  two equal rows. 
  Then either $\sat(n,F)$ is constant for large $n$, or $\sat(n,F)\ge
  n+1$ for every $n$.
\end{theorem}

\begin{proof}
  If some $M\in\Sat(n,F)$ has at most $n$ columns, then a well-known
  theorem of Bondy~\cite{bondy:72a} (see, \di{e.g.}, Theorem~2.1
  in~\cite{bollobas:c}) implies that there is $i\in[n]$ such that the
  removal of the $i$th row does not create two equal columns.  Since $F$
  has no two equal rows, the duplication of any row cannot create a
  forbidden submatrix, so $\sat(n+1,F)\ge\sat(n, F)$.  However, by the
  remark made just prior to the theorem, the duplication of the $i$th row
  gives an $(n+1)$-row $F$-saturated matrix, implying
  $\sat(n+1,F)\le\sat(n, F)$, as required.
\end{proof}

Suppose that $\C F$ consists of $k$-row matrices. Is there any good
general upper bound on $\forb(n,\C F)$ or $\sat(n,\C F)$? There were
different papers dealing with general upper bounds on $\forb(n,\C F)$,
for example, by Anstee and F\"uredi~\cite{anstee+furedi:86}, by Frankl,
F\"uredi and Pach~\cite{frankl+furedi+pach:87} and
by Anstee~\cite{anstee:95}, until the conjecture of Anstee and
F\"uredi~\cite{anstee+furedi:86} that $\forb(n,\C F)=O(n^k)$ for any
fixed $\C F$ was elegantly proved by F\"uredi
(see~\cite{anstee+griggs+sali:97} for a proof).

On the other hand, we can show that $\sat(n,\C F)=O(n^{k-1})$ for any
family $\C F$ of $k$-row matrices (including infinite families). Note that
the exponent $k-1$ cannot be decreased in general since, for example,
$\sat(n,T_k^k)=f(n,k-1)$.

\begin{theorem}\label{th::hsf:k-1}
For any family $\C F$ of $k$-row matrices,
$\sat(n,\C F)=O(n^{k-1})$.
\end{theorem}

\begin{proof}
  We may assume that $K_k$ is $\C F$-admissible (\textit{i.e.} every matrix
  of $\C F$ contains a pair of equal columns) for otherwise we are home
  by~\eqref{eq::sauer} as then $\sat(n,\C F)\le\forb(n,K_k)=O(n^{k-1})$.

  Let us define some parameters $l$, $d$, and $m$ that depend on $\C F$.
  Let $l=l(\C F)\in[0,k]$ be the smallest number such that there exists $s$
  for which $[sT_k^{\le l},T_k^{>l}]$ is not $\C F$-admissible. (Clearly,
  such $l$ exists: if we set $l=k$, then $sT_k^{\le l}=sK_k$ contains any
  given $k$-row submatrix for all large $s$.) Let $d=d(\C F)$ be the
  maximal integer such that $[sT_k^{<l},dT_k^l,T_k^{>l}]$ is $\C
  F$-admissible for every $s$. Note that $d\ge1$ as
  $[sT_k^{<l},T_k^l,T_k^{>l}]= [sT_k^{<l},T_k^{\ge l}]$ cannot contain a
  forbidden submatrix by the choice of $l$. Choose the minimal $m=m(\C
  F)\ge 0$ such that $[mT_k^{<l},(d+1)T_k^l,T_k^{>l}]$ is not $\C
  F$-admissible. The subsequent argument will be valid provided $n$ is
  large enough, which we shall tacitly assume.

  We consider the two possibilities $l(\C F)<k$ and $l(\C
  F)=k$ separately. Suppose first that $l(\C F)<k$. Consider the following set
  system:$$ H=\bigcup_{j\in[d-1]}\textstyle \left\{Y\in \subs{[n]}{l+1} :
  \sum_{y\in Y} y\equiv j\!\!\!\!\pmod{n}\right\}.$$ Here $\subs{X}{i}=\{Y\subseteq X :
  |Y|=i\}$ denotes the set of all subsets of $X$ of size $i$.

Note that any $A\in\subs{[n]}{l}$ is contained in at most $d-1$ members of $H$,
as there are at most $d-1$ possibilities to choose $i\in[n]\setminus
A$ so that $A\cup\{i\}\in H$: namely, $i\equiv j-\sum_{a\in A} a\pmod{n}$ for
$j\in[d-1]$.

On the other hand, the collection $H'$, of all $l$-subsets of $[n]$
contained in fewer than $d-1$ members of $H$, has size at most
$2(d-1)\binom{n}{l-1}$.  Indeed, if $A\in H'$ then, using the previous
observation, it must be that for some $j\in[d-1]$
and $x\in A$ we have $2x\equiv j-\sum_{a\in A\setminus\{x\}}a\pmod n$: hence, once
$A\setminus\{x\}$ and $j$ have been chosen, there are at most $2$ choices
for $x$.

Call $X\in\subs{[n]}{k}$ {\dd bad} if, for some $A\in \subs{X}{l}$,
\begin{equation*}\label{eq::nAB}
 |\{Y\in H : Y\cap X=A\}|\le
d-2.
\end{equation*}
To obtain a bad $k$-set $X$, we either complete some $A\in H'$ to any
$k$-set, or we take any $l$-set $A$ and let $X$ contain some member of $H$
that contains $A$. Therefore, the number of bad sets is at most$$
 2(d-1)\binom{n}{l-1}\binom{n}{k-l}+\binom{n}{l}(d-1)\binom{n}{k-l-1}=
O(n^{k-1}).$$

Let $M'=[N,T_n^l]$, where $N$ is the $n\times|H|$ incidence 
matrix of $H$. Then we have that 
 $$
 M'(X,)\subseteq[\textstyle e(M') T_k^{<l},dT_k^l,T_k^{l+1}],\quad
\mbox{for any $X\in\subs{[n]}{k}$.}$$
 Hence, $M'$ cannot contain a
forbidden submatrix by the definition of $d$. Now complete it to arbitrary
$M=[M',M'']\in\Sat(n,\C F)$ by adding new columns as long as no 
forbidden submatrix is created.

Suppose that $e(M'')\neq O(n^{k-1})$. Then, by~\eqref{eq::sauer}, $K_k\cong
M''(X,Y)$ for some $X,Y$. Now, remove the columns corresponding to $Y$
from $M''$ and repeat the procedure as long as possible to obtain
more than $O(n^{k-1})$ column-disjoint copies of $K_k$ in $M''$. No
$X\in\subs{[n]}{k}$ can appear more than $d$ times: otherwise (because
$T_n^l(X,)\supseteq mT_k^{<l}$ for all large $n$)
we have that $M(X,)=[M',M''](X,)\supseteq [mT_k^{<l},(d+1)K_k]$
is not $\C F$-admissible. Since we have
$O(n^{k-1})$ bad $k$-sets of rows and, by above, each has
at most $d$ column-disjoint copies of $K_k$, we have that $K_k\subseteq M''(X,)$
for at least
one \emph{good} (\emph{i.e.,} not bad) 
$X\in\subs{[n]}{k}$. But then $N(X,)\supseteq (d-1)T_k^l$. Moreover, since $T_n^l(X,)\supseteq [mT_k^{<l}, T_k^l]$ for all large $n$, we obtain
$$
M(X,) = [N,T_n^l, M''](X,)\supseteq [(d-1)T_k^l, mT_k^{<l}, T_k^l, K_k] = [(m+1)T_k^{<l}, (d+1)T_k^l, T_k^{>l}].
$$
Thus, $M(X,)$ contains a forbidden matrix.
This contradiction proves the required bound for $l<k$.

Consider now the other possibility, that $l=l(\C F)$ equals $k$. The above
argument does not work in this case because the size of $M'\supseteq T_n^l$
is too large. Let $\C F^*$ consist of those $k$-row matrices $F$ such that
$[dT_k^k,F]$ is not $\C F$-admissible, where $d=d(\C F)$.  Note that
$[sT_k^{<k},T_k^k]\in\C F^*$ for all large $s$ by the definition of
$d$. Thus $l(\C F^*)<k$ and by the above argument we can find
$L\in\Sat(n-d,\C F^*)$ with $O(n^{k-1})$ columns. Define
 $$
 M'=
 \left[
 \begin{array}{cc}
 dT_{n-d}^{n-d}&L\\
 T_d^1 & e(L)T^0_d
 \end{array}\right],
 $$
 that is, $M'$ is obtained from $[dT_{n-d}^{n-d},L]$ by adding 
$d$ extra rows that encode the sets $\{i\}$, $i\in [d]$. Note that
$M'$ does not have multiple columns even if $T_{n-d}^{n-d}$ is a column of~$L$ because $d\ge 1$.

Take arbitrary $X\in {[n]\choose k}$.  If $X\subseteq [n-d]$, then
$M'(X,)=[dT_k^k,L(X,)]$ is $\C F$-admissible because $L$ is $\C
F^*$-admissible; otherwise $M'(X,)\subseteq [e(M')T_k^{<k},T^k_k]$
is $\C F$-admissible because $l(\C F)=k$. Thus $M'$ is $\C F$-free.

Complete $M'$ to an arbitrary $M\in\Sat(n,\C F)$. Let $C$ be any added
column. Since 
 $$
 [M',C]([n-d],)=[dT_{n-d}^{n-d},L,C([n-d])]
 $$ 
 is $\C F$-free,
we have that $[L,C([n-d])]$ is $\C F^*$-free. By the $\C F^*$-saturation
of $L$, we have that $C([n-d])$
is a column of $L$.  Hence
 $$
 \sat(n,\C F)\le e(M)\le 2^d \, e(L)+d=O(n^{k-1}),
 $$
 proving the theorem.
\end{proof}

\begin{remark}  Theorem~\ref{th::hsf:k-1} is the
matrix analog of the main result in~\cite{pikhurko:99:cpc} that
$\sat(n,\C F)=O(n^{k-1})$ for any finite family $\C F$ of $k$-graphs.\end{remark}

\section{Forbidding Complete Matrices}\label{hsfcases}

Let us investigate the value of $\sat(n,K_k)$. (Recall that $K_k$ is the
$k\times 2^k$-matrix consisting of all distinct $k$-columns.) We are able
to settle the cases $k=2$ and $k=3$.

We will use the following trivial lemma a couple of times.

\begin{lemma}\label{lm::10}
Each row of any $M\in\Sat(n,K_k)$, $n\ge k$, contains at least
$l$ ones and at least $l$ zeros, where $l=2^{k-1}-1$.
\end{lemma}

\begin{proof}
Suppose on the contrary that the first row $M(1,)$ has $m_0$ zeros followed by
$m_1$ ones with $m_0\ge m_1$ and $l>m_1$.

For $i\in[m_0]$, let $C_i$ equal the $i$th column of $M$ with the first
entry $0$ replaced by $1$. Then the addition of $C_i$ to $M$ cannot create
a new copy of $K_k$, because the first row of $M'$ contains too few $1$'s,
while $C_i([2,n])$ is already a column of $M([2,n],)$, which does not
contain $K_k$.  Therefore, $C_i$ must be a column of $M$. Since $i\in
[m_0]$ was arbitrary, we have $m_0=m_1$.

But then $M$ has at most $2^k-2$ columns, which is a contradiction.
\end{proof}


\begin{theorem}\label{th::forb:k=2}
For $n\ge1$, we have $\sat(n,K_2)=n+1$.
\end{theorem}

\begin{proof} The upper bound is given by $T_n^{\le 1}\in\Sat(n,K_2)$.

Suppose that the statement is not true, that is, there exists a
$K_2$-saturated matrix with its size not exceeding its
order. By Theorem~\ref{th::hsf:const}, $\sat(n,K_2)$ is eventually
constant so we can find an $n\times m$-matrix $M\in\Sat(n,K_2)$ having
two equal rows for some $n\in\mathbb N$.

As we are free to complement and permute rows, we may assume that, for some
$i\ge2$, $M(1,)=\dots=M(i,)$ while $M(j,)\neq M(1,)$ and
$M(j,)\neq 1-M(1,)$ for any $j\in[i+1,n]$. Note that $i<n$ as we do not
allow multiple columns in $M$ (and $m\ge e(K_2)-1=3$).

Let $j\in[i+1,n]$. By Lemma~\ref{lm::10}, the $j$th row $M(j,)$ contains
both $0$'s and $1$'s. By the definition of $i$, $M(j,)$ is not equal to
$M(1,)$ nor to $1-M(1,)$. It easily follows that there are $f_j,g_j\in[m]$
with $M(1,f_j)=M(1,g_j)$ and $M(j,f_j)\neq M(j,g_j)$.  Again by
Lemma~\ref{lm::10}, we can furthermore find $h_j\in[m]$ with
$M(1,h_j)=1-M(1,f_j)$. Let $b_j=M(j,h_j)$. By exchanging $f_j$ and $g_j$ if
necessary, we can assume that $M(j,g_j)=b_j$.

Now, as $M\in\Sat(n,K_2)$, the addition of the column$$
 C=(1,(0)^{i-1},b_{i+1},\dots,b_n)^T$$
 (which is not in $M$ because $C(1)\neq C(2)$) must create a new
$K_2$-submatrix, say in the $x$th and $y$th rows for some $1\le x<y\le n$.
Clearly, $\{x,y\}\nsubseteq [i]$ because each column of $M([i],)$ is
either $((0)^i)^T$ or $((1)^i)^T$. Also, it is impossible that $x\in[i]$
and $y\in[i+1,n]$ because then, for some $a_1,a_2\in[m]$,
$M(y,a_1)=M(y,a_2)=1-C(y)=1-b_y$, $M(x,a_1)=1-M(x,a_2)$ and we can see
that $K_2$ is isomorphic to  $M(\{x,y\},\{a_1,a_2,g_y,h_y\})$, which
contradicts $K_2\nsubseteq  M(\{x,y\},)$. So we have to assume that
$i<x<y\le n$. 

As $K_2\nsubseteq M(\{x,y\},)$, no column of $M(\{x,y\},)$ can equal
$C(\{x,y\})=(b_x,b_y)^T$.  In particular, since $M(x,g_x)=M(x,h_x)=b_x$ and
similarly for~$y$, we must have $\{g_x,h_x\}\cap \{g_y,h_y\}=\emptyset$,
and moreover $M(y,g_x)=M(y,h_x)=1-b_y$. But then
\begin{equation*}\label{eq::hsfcases:ghgh}
 K_2\cong M(\{1,y\},\{g_x,h_x,g_y,h_y\}),
\end{equation*}
 which is a contradiction proving our theorem.
\end{proof}

Note that $\forb(n,K_2)=n+1$ for $n\ge 1$; the upper bound follows, for
example, from Formula~\eqref{eq::sauer} with $k=2$.  Thus
Theorem~\ref{th::forb:k=2} yields that $\sat(n,K_2)=\forb(n,K_2)$ which, in
our opinion, is rather surprising. A greater surprise is yet to come as we
are going to show now that $\sat(n,K_3)$ is constant for $n\ge4$.

\begin{theorem}\label{th::K3}
For $K_3$ the following holds:
$$
\sat(n,K_3)=
\begin{cases}
7, &\text{if $n=3$,}\\
10, &\text{if $n\ge 4$.}
\end{cases}
$$
\end{theorem}

\begin{proof}
 The claim is trivial for $n=3$, so assume $n\ge 4$.
A computer search~\cite{dudek+pikhurko+thomason:09:arxiv} revealed that
\begin{equation*}\label{eq::agt}
 \sat(4,K_3)=\sat(5,K_3)=\sat(6,K_3)=\sat(7,K_3)=10,
\end{equation*}
 which suggested that $\sat(n,K_3)$ is constant. An example of a
$K_3$-saturated $6\times 10$-matrix is the following.$$
M=\left[\begin{array}{llllllllll}
 0 & 0 & 0 & 0 & 1 & 1 & 0 & 1 & 1 & 1\\
 0 & 0 & 1 & 1 & 0 & 0 & 0 & 1 & 1 & 1\\
 0 & 1 & 0 & 1 & 0 & 0 & 1 & 0 & 1 & 1\\
 1 & 0 & 0 & 0 & 0 & 1 & 1 & 0 & 1 & 1\\
 1 & 0 & 1 & 0 & 0 & 0 & 1 & 1 & 0 & 1\\
 0 & 1 & 0 & 0 & 1 & 0 & 1 & 1 & 0 & 1\end{array}\right]. $$
 It is possible (but very boring) to check by hand that $M$ is indeed
$K_3$-saturated as is, in fact, any $n\times 10$-matrix $M'$ obtained
from $M$ by duplicating any row, \di{cf.}\
Theorem~\ref{th::hsf:const}. (The symmetries of $M$ shorten the
verification.)
A $K_3$-saturated $5\times10$-matrix can be obtained from $M$ by
deleting one row (any). For $n=4$, we have to provide a special
example:$$
 M=\left[\begin{array}{llllllllll}
 0 & 0 & 0 & 0 & 0 & 0 & 0 & 1 & 1 & 1\\
 0 & 0 & 0 & 0 & 1 & 1 & 1 & 0 & 1 & 1\\
 0 & 0 & 1 & 1 & 0 & 0 & 1 & 1 & 0 & 1\\
 0 & 1 & 0 & 1 & 0 & 1 & 0 & 1 & 1 & 0\end{array}\right]. $$

So $\sat(n,K_3)\le 10$ for each $n\ge 4$ and, to prove the theorem, we
have to show that no $K_3$-saturated matrix $M$ with at most $9$ columns
and at least $4$ rows can exist. Let us assume the contrary.

\claimskip
\noindent{\em Claim~1} Any row of $M\in\Sat(n,K_3)$ necessarily contains
at least four 0's and at least four 1's, for $n\ge4$.\claimskip

\bcpf Suppose, contrary to the claim, that the first row $M(1,)$
contains only three 0's, say in the first three columns. (By
Lemma~\ref{lm::10} we must have at least three 0's.)

If we replace the $i$th of these 0's by 1, $i\in[3]$, then the obtained
column $C_i$, if added to $M$, does not create any $K_3$-submatrix. Indeed,
the first row of $[M,C_i]$ contains at most three 0's, while $C_i([2,n])$
is a column of $M([2,n],)\not\supseteq K_3$.  As $M$ is $K_3$-saturated,
$C_1$, $C_2$ and $C_3$ are columns of $M$.  These columns differ only in
the first entry from $M(,1)$, $M(,2)$ and $M(,3)$ respectively. Therefore,
for each $A\in\subs{[2,n]}{3}$, the matrix $M(A,)$ can contain at most
$e(M)-3\le 6$ distinct columns. But then any column $C$ which is not a
column of $M$ and has top entry $1$ ($C$ exists as $n\ge4$) can be added
to~$M$ without creating a $K_3$ submatrix, because the first row of $[M,C]$
contains at most three 0's. This contradiction proves Claim~1.\ecpf

Therefore, $e(M)$ is either 8 or 9. As we are free to complement the rows,
we may assume that each row of $M$ contains exactly four $1$'s. Call
$A\in\subs{[n]}{3}$ (and also $M(A,)$) {\dd nearly complete} if $M(A,)$ has
7 distinct columns.

\claimskip\noindent{\em Claim~2} Any nearly complete $M(A,)$ contains
$(0,0,0)^T$ as a column.\claimskip

\bcpf Indeed, otherwise $M(A,)\supseteq T_3^{\ge1}$ which already contains four
$1$'s in each row; this implies that the (one or two) remaining columns
must contain zeros only. Hence $M(A,)\supseteq K_3$,  
which is a contradiction.\ecpf

\noindent{\em Claim~3} Every nearly complete $M(A,)$ contains $T_3^1$ as a
submatrix.\claimskip

\bcpf Indeed, if $(0,0,1)^T$ is the missing column of $M(A,)$, then some 7
columns contain a copy of $K_3\setminus(0,0,1)^T$. By counting 1's in the
rows we deduce that the remaining column(s) of $M(A,)$ must have exactly
one non-zero entry, and moreover one of these columns equals $(0,0,1)^T$,
which is a contradiction.\ecpf

By the $K_3$-saturation of $M$ there exists some nearly complete $M(A,)$;
choose one such. Assume without loss of generality that $A=[3]$ and
that the first 7 columns of $M([3],)$ are distinct. We know that the
$3$-column missing from $M([3],[7])$ has at least two $1$'s.

If $(1,1,1)^T$ is missing, then $M([3],[7])$ contains exactly three
ones in each row, so the remaining column(s) of $M$ must contain an
extra 1 in each row. As $(1,1,1)^T$ is the missing column, we conclude
that $e(M)=9$ and the 8th and 9th columns of $M([3],)$ are, up to a
row permutation, $(0,0,1)^T$ and $(1,1,0)^T$. This implies that
$M([3],)$ contains the column $(0,0,0)^T$ only once. Thus at least one
of the columns $C_0=((0)^n)^T$ and $C_1=((0)^{n-1},1)^T$ is not in
$M$ and its addition creates a copy of $K_3$, say on the rows indexed
by $B\in \subs{[n]}3$. The submatrix $M(B,)$ is nearly complete
and,  by Claims~2 and~3, contains $T_3^{\le 1}$. But
both $C_0(B)$ and $C_1(B)$ are columns of $T_3^{\le 1}
\subseteq M(B,)$, which is a contradiction.

Similarly, if $(1,1,0)^T$ is missing, then one can deduce that $e(M)=9$ and, 
up to a row
permutation, $M([3],\{8,9\})$ consists of the columns
$(1,0,0)^T$ and $(0,1,0)^T$. Again, the column $(0,0,0)^T$
appears only once in $M([3],)$, which leads to a contradiction as above, 
completing the proof of the theorem.
\end{proof}

We do not have any non-trivial results concerning $K_k$, $k\ge4$, except
that a computer search~\cite{dudek+pikhurko+thomason:09:arxiv} showed that
$\sat(5,K_4)=22$ and $\sat(6,K_4)\le 24$. (We do not know if a
$K_4$-saturated $6\times 24$-matrix discovered by a partial search is
minimum.)

\begin{problem}
For which $k\ge 4$, is $\sat(n,K_k)=O(1)$?
\end{problem}

\section{Forbidding Small Matrices}\label{smalls}

In this final section we try to gain further insight into the
$\sat$-function by computing $\sat(n,F)$ for some forbidden matrices with
up to three rows.
\subsection{Forbidding $1$-Row Matrices}

For any given $1$-row matrix $F$, we can determine $\sat(n,F)$ for all but
finitely many values of $n$. The answer is unpleasantly intricate.

\begin{proposition}\label{pr::1-row} Let $F=((0)^m,(1)^l)=[mT_1^0,lT_1^1]$ with $l\ge m$.
Then, for $n\ge \max(l-1,1)$,
$$
\sat(n,F)=
\begin{cases}
l,& \text{if $m=0$ and $l\le2$ or if $m=1$ and $l\ge 1$ is a power of $2$},\\
l+1,   &\text{if $m=0$ and $l\ge 3$ or if $m=1$ and $l$ is not a power of $2$,}\\
l+m-1, &\text{if $m\ge 2$ and $l\ge 2$.}
\end{cases}
$$
 \end{proposition}

\begin{proof} Assume that $l\ge3$, as the case $l\le 2$ is trivial.

For $m\in\{0,1\}$ an example of $M\in\Sat(n,F)$ with $e(M)=l+1$ can be built 
by taking $T_n^0$, $T_n^n$, $\chi_{[l-2]}$, and $\chi_{[n]\setminus \{i\}}$ for
$i\in[l-2]$ as the columns. If $m=1$ and $l=2^k$, one can do slightly 
better by adding $n-k$ copies of the row $((1)^{l})$ to $K_k$.

Let us prove the lower bound for $m\in\{0,1\}$. Suppose that some
$F$-saturated matrix $M$ has $n\ge l-1$ rows and $c\le l$ columns. First,
let $m=0$.  As $c<2^n$ and $M$ contains the all-$0$ column, we have $c=l$
and some row $M(i,)$ contains exactly $l-1$ ones. As we are not allowed
multiple columns in $M$, some other row, say $M(j,)$, has at most $l-2$
ones. Then $\chi_{\{j\}}$ is not a column of $M$ but its addition does not
create $l$ ones in a row, a contradiction. Let $m=1$.  Trivially, $e(M)\ge
e(F)-1=l$. It remains to show that $l$ is a power of $2$ if $e(M)=l$.  Let
$C$ be the column whose $i$th entry is $1$ if and only if
$M(i,)=(1)^l$. Then the addition of the column $C$ cannot create an
$F$-submatrix, and so $C$ is already a column of $M$. Let
$C=M(,1)=((0)^i,(1)^{n-i})^T$. The last $n-i$ rows of $M$ consist of $1$'s
only. Since $l\ge 3$ and $M$ has no multiple columns, we have that $i\ge 2$
and that $M([i],[2,l])$ must contain at least one~$0$, say
$M(i,l)=0$. Since the addition of $\chi_{[i,n]}$ cannot create $F$, it is
already a column of $M$. Thus each row of $M([i],)$ has at least two 0's,
and to avoid a contradiction we must have $M([i],)\cong K_i$ and $l=2^i$.
This completes the case when $m\le 1$.

For $m\ge 2$, let $M$ consist of $T_n^n$ plus $\chi_{\{i\}}$, $i\in[m-2]$,
plus $\chi_{[n]\setminus \{i\}}$, $i\in[l-1]$ and $\chi_{[m-1,l-1]}$.
Clearly, each row of $M$ contains $l$ 1's and $m-1$ 0's, so the addition
of any new column (which must contain at least one $0$) creates an
$F$-submatrix and the upper bound follows. The lower bound is trivial.
\end{proof}

\begin{remark} The case when $n\le l-2$ in\ Proposition~\ref{pr::1-row} seems
  messy so we do not investigate it here.\end{remark}

\subsection{Forbidding $2$-Row Matrices}
Now let us consider some particular $2$-row matrices.

Let $F=lT_2^2$, that is, $F$ consists of the column $(1,1)^T$ taken $l$
times. Trivially, for $l=1$ or $2$, $\sat(n,lT_2^2)=n+l$, with $T_n^{\le
  1}$ and $[T_n^{\le 1},T_n^n]$ being the only extremal matrices.  For
$l\ge3$, we can only show the following lower bound. It is almost sharp
for $l=3$, when we can build a $3T_2^2$-saturated $n\times (2n+2)$-matrix by
taking $T_n^{\le1}$, $\chi_{[n-1]}$, $\chi_{[n]}$, plus $\chi_{\{i,n\}}$
for $i\in[n-1]$.

\begin{lemma}\label{lm::lT22} For $l\ge 3$ and $n\ge 3$, $\sat(n,lT_2^2)\ge 2n+1$.
\end{lemma}

\begin{proof}
  Let $M=[T_n^{\le1},M']$ be $lK_2^2$-saturated. Note that $M'$ must have
  the property that every column $\chi_A$, with $A\in\subs{[n]}{2}$, either
  belongs already to $M'$, or its addition creates an $F$-submatrix; in the
  latter case, exactly $l-1$ columns of $M'$ have ones in both positions of
  $A$.  Therefore, by adding to $M'$ some columns of $T_n^2$ (with possibly
  some columns being added more than once), we can obtain a new matrix
  $M''$ such that, for every $A\in\subs{[n]}{2}$, $M''(A,)$ contains the
  column $(1,1)^T$ exactly $l-1$ times. If we let the set $X_i$ be encoded
  by the $i$th row of $M''$ as its characteristic vector, we have that
  $|X_i\cap X_j|=l-1$ for every $1\le i<j\le n$.  The result of
  Bose~\cite{bose:49} (see \cite[Theorem 14.6]{jukna:ec}), which can be
  viewed as an extension of the famous Fisher inequality~\cite{fisher:40},
  asserts that, either the rows of $M''$ are linearly independent over the
  reals, or $M''$ has two equal rows, say $X_i=X_j$.  The second case is
  impossible here, because then $|X_i|=l-1$ and each other $X_h$ contains
  $X_i$ as a subset; this in turn implies that the column $((1)^n)^T$
  appears at least $l-1\ge 2$ times in $M''$ and (since $n\ge 3$) the same
  number of times in $M'$, a contradiction.  Thus the rank of $M''$ over
  the reals is $n$. Note that every column $C\in T_n^2$ added to $M'$
  during the construction of $M''$ was already present in $M'$ (otherwise
  $C$ contradicts the assumption that $M$ is $lT_2^2$-saturated). Thus the
  matrices $M'$ and $M''$ have the same rank over the reals. We conclude
  that $M'$ has at least $n$ columns and the lemma follows.
\end{proof}

Let us show that Lemma~\ref{lm::lT22} is sharp for $l=3$ and 
some $n$. Suppose there
exists a \emph{symmetric $(n,k,2)$-design} (meaning we have $n$ $k$-sets
$X_1,\dots,X_n\in \subs{[n]}{k}$ such that every pair $\{i,j\}\in 
\subs{[n]}{2}$ is
covered by exactly two $X_i$'s). Let $M$ be the $n\times n$-matrix whose rows
are the characteristic vectors of the sets $X_i$. Then
$[T_n^{\le 1}, M]$ is a $3T_2^2$-saturated $n\times(2n+1)$-matrix. For
$n=4$, we can take all $3$-subsets of $[n]$. For $n=7$, we can take the family
$\{[7]\setminus Y_i : i\in [7]\}$, where $Y_1,\dots,Y_7\in \subs{[7]}{3}$ form
the Fano plane. Constructions of such designs for $n=16$, $37$, $56$, and $79$
can be found in \cite[Table 6.47]{colbourn+dinitz:hcd}.

Of course, the non-existence of a symmetric $(n,k,2)$-design does not
directly imply anything about $\sat(n,3T_2^2)$, since a minimum
$3T_2^2$-saturated matrix $[T_n^{\le 1}, M]$ need not have the same number
of ones in the rows of $M$.

Lemma~\ref{lm::lT22} is not always optimal for $l=3$. One trivial example is
$n=3$. Another one is $n=5$.

\begin{lemma}\label{lm::3T22} $\sat(5,3T_2^2)=12$.\end{lemma}

\begin{proof} Suppose, on the contrary, that we have a $3T_2^2$-saturated
  $5\times (s+6)$-matrix $M=[N,T_5^{\le 1}]$ with $s\le 5$.  Let
  $X_1,\dots,X_5$ be the subsets of $[s]$ encoded by the rows of $N$.

  If, for example, $X_1=[s]$, then every $X_i$ with $i\ge 2$ has at most
  two elements. Let $C_1=(0,1,1,0,0)^T$, $C_2=(0,0,0,1,1)^T$ and
  $C_3=(0,0,1,1,0)^T$. None of these columns is in $M$ so the addition of
  any one of them creates a copy $3T_2^2$. So we may assume that
  $M(\{2,3\},\{a,b\})=M(\{4,5\},\{c,d\})=M(\{3,4\},\{e,f\})=2T_2^2$. If
  $\{a,b\}=\{c,d\}$ then $M(,a)$ and $M(,b)$ are two equal columns with
  all~1's, a contradiction. Hence $\{a,b\}\ne\{c,d\}$, and so at least one
  of $\{e,f\}\ne \{a,b\}$ or $\{e,f\}\ne \{c,d\}$ holds: we may assume the
  former. But then $M(\{1,3\},)$ contains $3T_2^2$, a contradiction.

  Thus we can assume that each $X_i$ with $i\in [5]$ has at most $s-1$
  elements. If $X_1\subseteq \{1,2\}$, then by considering columns that
  begin with $1$ and have one other entry $1$, we conclude that
  $X_1=\{1,2\}$ and that every $X_i$ contains $X_1$ as a subset. Thus
  $M(,\{1,2\})=2T_5^5$, that is, $M$ has two equal columns, a
  contradiction.

  So we can assume that each $|X_i|\ge 3$, which also implies that $s=5$.
  If $X_1=[4]$, then for each $i\in [2,5]$ we have $5\in X_i$ (because
  $|X_i|\ge 3$ and $M$ is $3T_2^2$-free). Each two of the sets
  $X_2,\dots,X_5$ have to intersect in exactly two elements, which is
  impossible.

  Thus each $|X_i|=3$. A simple case analysis gives a contradiction in this
  case as well.\end{proof}

\begin{problem} Determine $\sat(n,3T_2^2)$ for every $n$.\end{problem} 


\begin{remark}
  It is interesting to note that if we let $F=[lT_2^2,(0,1)^T]$ then
  $\sat(n,F)$-function is bounded.  Indeed, complete
  $M'=[\chi_{[n]\setminus\{i\}}]_{i\in[l]}$ to an arbitrary $F$-saturated
  matrix $M$. Clearly, in any added column all entries after the $l$th
  position are either 0's or 1's; hence $\sat(n,F)\le 2\cdot 2^l$.
\end{remark}

It is easy to compute $\sat(n,T_2^1)$ by observing that the $n$-row matrix
$M_Y$ whose columns encode $Y\subseteq2^{[n]}$ is $T_2^1$-free if and only
if $Y$ is a chain --- that is, for any two members of $Y$, one is a subset
of the other. Thus $M_Y$ is $T_2^1$-saturated if and only if $Y$ is a
maximal chain without repeated entries.  As all maximal chains in $2^{[n]}$
have size $n+1$, we conclude that$$ \sat(n,T_2^1)=\forb(n,T_2^1)=n+1,\quad
n\ge2.$$

\begin{theorem}\label{th::sat=3} Let $F=[T_2^0,T_2^2]=\left[\begin{array}{ll}
 0 & 1\\ 0 & 1\end{array}\right]$. Then $\sat(n,F)=3$, $n\ge2$.
\end{theorem}

\begin{proof}
For $n\ge3$, the matrix $M$ consisting of the columns
$(0,1,(1)^{n-2})^T$, $(1,0,(1)^{n-2})^T$ and $(0,0,(1)^{n-2})^T$ can be
easily verified to be $F$-saturated and the upper bound follows.

Since $n=2$ is trivial, let $n\ge 3$.  Any $2$-column $F$-free matrix $M$
is, without loss of generality, the following: we have $n_{00}$ rows
$(0,0)$, followed by $n_{11}$ rows $(1,1)$, $n_{10}$ rows $(1,0)$ and
$n_{01}$ rows $(0,1)$, where $n_{10}\le1$ and $n_{01}\le1$. Since (by
taking complements if necessary) we may assume $n_{00}\le n_{11}$, we have
$n_{11}\ge1$ because $n\ge3$. But then the addition of a new column
$((0)^{n_{00}+1},1,1,\dots)^T$ does not create an $F$-submatrix.
\end{proof}

\begin{theorem}\label{th::T2ge1} Let $F=T_2^{\ge1}=\left[\begin{array}{lll}
 0 & 1 &1\\ 1 & 0 & 1\end{array}\right]$. Then$$
 \sat(n,F)=\forb(n,F)=n+1,\quad n\ge 2.$$
\end{theorem}

\begin{proof}
Clearly, $\forb(n,F)\le\forb(n,K_2) =n+1$.

Suppose the theorem is false and that $\sat(n,F)\le n$ for some~$n$. Since
the rows of $F$ are distinct, Theorem~\ref{th::hsf:const} shows that 
$\sat(n,F)$ is bounded.

It follows that, if $n$ is large enough, then $M\in\Sat(n,F)$ has two equal
rows, for example, $M(1,)=M(2,)=((1)^l,(0)^m)$. By considering the column
$(1,0,\dots,0)^T$ that is not in $M$, we conclude that $l,m\ge 1$.  Let
$X=[l]$ and $Y=[l+1,l+m]$. Define
$$
 A_i=\{j\in[l+m] : M(i,j)=1\},\quad  i\in[n].$$
 (For example, $A_1=A_2=X$.)
As $M$ is $F$-free, for every $i,j\in[n]$, the sets $A_i$ and $A_j$ are
either disjoint or one is a subset of the other. For $i\in[3,n]$, let
$b_i=1$ if $A_i$ strictly contains $X$ or $Y$ and let $b_i=0$ otherwise
(that is, when $A_i$ is contained in $X$ or $Y$). Let $b_1=1$ and $b_2=0$.

Clearly, $C=(b_1,\dots,b_n)^T$ is not a column of $M$ so its addition
creates a forbidden submatrix, say $F\subseteq[M,C](\{i,j\},)$. Of
course, $b_i=b_j=0$ is impossible because $(0,0)^T\nsubseteq  F$. If
$b_i=b_j=1$ then necessarily $A_i\cap A_j\neq \emptyset$, and
$M(\{i,j\},)\supseteq (1,1)^T$ contains $F$, a contradiction.  Finally, if
$b_i\neq b_j$, \di{e.g.}, $b_i=0$, $b_j=1$ and $i<j$, then $A_i\supseteq A_j$
(as $(0,1)^T$ cannot be a column of $M(\{i,j\},)$), which implies
$A_i=A_j$; but then we do not have a copy of $F$ as $(1,0)^T$ is
missing. This contradiction completes the proof.
\end{proof}

\begin{remark}
  It is trivial that
  $\sat(n,[(0,1)^T,(1,1)^T])=\sat(n,[(0,0)^T,(0,1)^T,(1,1)^T])=2$. We have
  thus determined the $\sat$-function for every simple $2$-row matrix.
\end{remark}

\subsection{Forbidding $3$-Row Matrices}

Here we consider some particular 3-row matrices. First we solve completely the
case when $F=[T_3^0,T_3^3]$. 

\begin{theorem} Let $F=[T_3^0,T_3^3]=
\left[\begin{array}{ll}
 0 & 1\\ 0 & 1\\ 0 & 1
 \end{array}
 \right]$. 
Then 
$$
\sat(n,F)=
\begin{cases}
7, &\text{if $n=3$ or $n\ge 6$,}\\
10, &\text{if $n=4$ or $5$.}
\end{cases}
$$
\end{theorem}

\begin{proof}
For the upper bound we define the following family of matrices.
$$
 M_4=\left[
 \begin{array}{llllllllll}
1 & 0 & 1 & 0 & 1 & 0 & 1 & 1 & 0 & 0\\
0 & 1 & 1 & 0 & 0 & 1 & 1 & 0 & 1 & 0\\
0 & 0 & 0 & 1 & 1 & 1 & 1 & 0 & 0 & 1\\
0 & 0 & 0 & 0 & 0 & 0 & 0 & 1 & 1 & 1
\end{array}
\right] 
$$

$$
 M_5=\left[
 \begin{array}{llllllllll}
1 & 1 & 0 & 1 & 1 & 0 & 1 & 0 & 1 & 0\\
1 & 0 & 1 & 1 & 0 & 1 & 0 & 1 & 1 & 0\\
0 & 1 & 1 & 1 & 0 & 0 & 1 & 1 & 0 & 1\\
0 & 0 & 0 & 0 & 1 & 1 & 1 & 1 & 0 & 0\\
0 & 0 & 0 & 0 & 0 & 0 & 0 & 0 & 1 & 1
\end{array}
\right]
$$

$$
 M_6=\left[
 \begin{array}{lllllll}
1 & 0 & 1 & 0 & 0 & 1 & 0\\
1 & 0 & 0 & 1 & 1 & 0 & 0\\
0 & 1 & 1 & 0 & 1 & 0 & 0\\
0 & 1 & 0 & 1 & 0 & 1 & 0\\
0 & 0 & 1 & 1 & 0 & 0 & 1\\
0 & 0 & 0 & 0 & 1 & 1 & 1
\end{array}
\right]
$$
For any $n\ge 7$ define the $(n\times 7)$-matrix $M_n$ by $M_n([6],) = M_6$ 
and
$M_n(i,) = \left[\begin{array}{lllllll}
0 & 0 & 0 & 0 & 0 & 0 & 0
 \end{array}\right]$ for every $7\le i\le n$.
 A computer search~\cite{dudek+pikhurko+thomason:09:arxiv} showed that $M_n$ is
a minimum $F$-saturated matrix for $3\le n \le 10$. This implies that
each $M_n$ with $n\ge 11$ is $F$-saturated. 
It remains to show that 
\begin{equation*}\label{eq:claim:zerorow}
\sat(n,F)\ge 7
\end{equation*}
for $n\ge 11$. In order to see this, we show the following result first.

\claimskip\noindent{\em Claim} If $M$ is an $F$-saturated  $n\times m$-matrix
with $n\ge 11$ and $m\le 6$ then $M$ contains a row with all zero entries or
with all one entries.\medskip

\bcpf Suppose, on the contrary, that we have a counterexample $M$.  We may
assume that the first 6 entries of the first column of $M$ are equal
to~0. Consider a matrix $A=M([6], \{2,\dots,m\})$. Note that every column
of $A$ contains at most two entries equal to 1, otherwise $M([6],)\supseteq
F$. Hence, the number of 1's in $A$ is at most $2(m-1)$. By our assumption,
each row of $A$ has at least one~$1$. Since $2(m-1)<12$, $A$ has a row with
precisely one~1. We may assume that $A(1,1)=1$ and
$A(1,i)=0$ for $2\le i\le m-1$. Let~$C_2$ be the second column of~$M$
(remember that $C_2(1)=A(1,1)=1$).

Consider the $n$-column $C_3=[0, C_2(\{2,\dots,n\})^T]^T$ which is obtained from
$C_2$ by changing the first entry to $0$. If it is not in $M$, then
$F\subseteq [M,C_3]$. This copy of $F$ has to contain the entry  
in which $C_3$ differs from $C_2$. But the only non-zero entry in Row 1 is
$M(1,2)$; thus $F\subseteq [C_2,C_3]$, which is an obvious contradiction. Thus
we may assume that $C_3$ is the third column of~$M$.

We have to consider two cases. First, suppose that $C_2(\{2,\dots,6\})$ has at
least one entry equal to $1$. Without loss of generality, assume that
$C_2(2)=C_3(2)=1$.

It follows that  $C_2(i)=C_3(i)=0$ for
$3\le i \le 6$ (otherwise the first and the second columns of $M$ would contain
$F$). Let 
 \begin{equation}\label{B}
 B=M(\{3,4,5,6\},\{4,\dots,m\}).
 \end{equation}
 By our assumption, each row of $B$ has at least one $1$; in particular $m\ge
5$.  Clearly,  $B$ contains at most $2(m-3)<8$ ones. Thus, by permuting
Rows~$3,\dots,6$ and Columns~$4,\dots,m$, we can assume that $B(1,1)=1$ while
$B(1,i)=0$ for $2\le i\le m-3$. 
Let $C_4$ be the fourth column of~$M$ and $C_5$ be
such that $C_4$ and $C_5$ differ at the third position only, \textit{i.e.},
$C_4(3)=1$ and $C_5(3)=0$. As before, $C_5$ must be in $M$, say it is the fifth
column.  Since
$C_4(\{4,5,6\})$ has at most one 1, assume that
$C_4(5)=C_4(6)=C_5(5)=C_5(6)=0$. We need another column $C_6$ with
$C_6(5)=C_6(6)=1$ (otherwise the fifth or the sixth row of $M$ would consist
of all zero entries). In particular, $m=6$.
But now the new column $C_7$ which differs from~$C_6$ at the fifth position
only (\textit{i.e.} $C_7(5)=0$ and $C_7(i)=C_6(i)$ for $i\neq 5$) should be
also in $M$, since $M$ is $F$-saturated. This contradicts $e(M)=6$. Thus the first case does not hold.

In the second case, we have $C_2(i)=C_3(i)=0$ for every $2\le i\le 6$. We may
define $B$ as in (\ref{B}) and get a contradiction in the same way as above.
This proves the claim.\ecpf

Suppose, contrary to the theorem, that we can find an $F$-saturated matrix
$M$ with $n\ge 11$ rows and $m\le 6$ columns. By the claim, $M$ has a
constant row; we may assume that the final row of $M$ is all zero, and let
$N=M([n-1],)$. If $C$ is an $(n-1)$-column missing from $N$, then the
column $Q=(C^T,0)^T$ is missing in $M$. Moreover, a copy of $F$ in $[M,Q]$
cannot use the $n$-th row. Thus $F\subseteq [N,C]$, which means that $N\in
\Sat(n-1,F)$ and $\sat(n-1,F)\le m\le 6$.  Repeating this argument, we
eventually conclude that $\sat(10,F)\le 6$, a contradiction to the results
of our computer search. The theorem is proved.\end{proof}

\begin{theorem} Let $F=[T_3^0,T_3^2,T_3^3]=
\left[\begin{array}{lllll}
0 & 0 & 1 & 1 & 1\\
0 & 1 & 0 & 1 & 1\\
0 & 1 & 1 & 0 & 1
 \end{array}
 \right]$.
Then
$$
\sat(n,F)=
\begin{cases}
7, &\text{if $n=3,6$ or $7$,}\\
9, &\text{if $n=4$ or $5$.}
\end{cases}
$$
Moreover, for any $n\ge 8$, $\sat(n,F)\le 7$.
\end{theorem}

\begin{proof} We define the following matrices:
$$
 M_4=\left[
 \begin{array}{lllllllll}
1 & 0 & 1 & 0 & 1 & 0 & 0 & 0 & 1\\
0 & 1 & 1 & 0 & 0 & 1 & 0 & 1 & 1\\
0 & 0 & 0 & 1 & 1 & 0 & 1 & 1 & 1\\
0 & 0 & 0 & 0 & 0 & 1 & 1 & 1 & 1 
\end{array}
\right],
$$

$$
 M_5=\left[
 \begin{array}{lllllllll}
1 & 1 & 1 & 0 & 1 & 0 & 1 & 0 & 1\\
0 & 1 & 0 & 1 & 0 & 1 & 0 & 1 & 1\\
0 & 0 & 1 & 0 & 1 & 0 & 1 & 1 & 1\\
0 & 0 & 0 & 1 & 1 & 0 & 0 & 1 & 1\\
0 & 0 & 0 & 0 & 0 & 1 & 1 & 1 & 1 
\end{array}
\right],
$$

$$
 M_6=\left[
 \begin{array}{lllllll}
1 & 1 & 0 & 0 & 1 & 1 & 0\\
1 & 0 & 1 & 1 & 0 & 1 & 0\\
1 & 0 & 1 & 0 & 1 & 0 & 1\\
0 & 1 & 1 & 1 & 1 & 0 & 0\\ 
0 & 1 & 1 & 0 & 0 & 1 & 1\\
0 & 0 & 0 & 1 & 1 & 1 & 1
\end{array}
\right],
$$
For any $n\ge 7$ let $M_n([6],) = M_6$
and
$M_n(i,) = \left[\begin{array}{lllllll}
0 & 0 & 0 & 1 & 1 & 1 & 1
 \end{array}\right]$ for every $7\le i\le n$ (\textit{i.e.} the last row of~$M_6$ is repeated $(n-6)$ times).
For $n=3,\dots,7$ the theorem (with $M_n$ being a minimum $F$-saturated matrix) follows from a computer search~\cite{dudek+pikhurko+thomason:09:arxiv}.
It remains to show that $M_n$, $n\ge 8$, is $F$-saturated. Clearly, this is the case, since $M_7$ is $F$-saturated and $F$ contains no pair of equal rows. 
\end{proof}

\begin{conjecture}
Let $F=[T_3^0,T_3^2,T_3^3]$. Then $\sat(n,F)=7$ for every $n\ge 8$.
\end{conjecture}

\begin{theorem} Let $F=T_3^{\le 2}=
\left[\begin{array}{lllllll}
0 & 1 & 0 & 0 & 0 & 1 & 1\\
0 & 0 & 1 & 0 & 1 & 0 & 1\\
0 & 0 & 0 & 1 & 1 & 1 & 0
 \end{array}
 \right]$.
Then
$$
\sat(n,F)=
\begin{cases}
7, &\text{if $n=3$,}\\
10, &\text{if $4\le n\le 6$.}
\end{cases}
$$
Moreover, for any $n\ge 7$, $\sat(n,F)\le 10$.
\end{theorem}

\begin{proof}
For $n=3,\dots,6$ the statement follows from a computer search~\cite{dudek+pikhurko+thomason:09:arxiv} with the following $F$-saturated matrices.
$$
 M_4=\left[
 \begin{array}{llllllllll}
0 & 1 & 0 & 1 & 0 & 1 & 1 & 0 & 0 & 1\\
0 & 0 & 1 & 1 & 0 & 0 & 1 & 1 & 0 & 1\\
0 & 0 & 0 & 0 & 1 & 1 & 1 & 0 & 1 & 1\\
0 & 0 & 0 & 0 & 0 & 0 & 0 & 1 & 1 & 1
\end{array}
\right]
$$

$$
 M_5=\left[
 \begin{array}{llllllllll}
1 & 0 & 1 & 0 & 1 & 0 & 0 & 0 & 1 & 1\\
0 & 1 & 0 & 1 & 1 & 0 & 1 & 0 & 0 & 1\\
0 & 0 & 1 & 0 & 1 & 1 & 1 & 0 & 0 & 1\\
0 & 0 & 0 & 1 & 1 & 0 & 0 & 1 & 1 & 1\\
0 & 0 & 0 & 0 & 0 & 1 & 1 & 1 & 1 & 1
\end{array}
\right]
$$
For any $n\ge 6$ let $M_n([5],) = M_5$
and
$M_n(i,) = \left[\begin{array}{llllllllll}
1 & 1 & 0 & 0 & 0 & 0 & 1 & 0 & 1 & 1
 \end{array}\right]$ for every $6\le i\le n$.
It remains to show that $M_n$, $n\ge 7$, is $F$-saturated. Clearly, this is the case, since $M_6$ is $F$-saturated and $F$ contains no pair of equal rows. 
\end{proof}

\begin{conjecture}
Let $F=T_3^{\le 2}$. Then $\sat(n,F)=10$ for every $n\ge 7$.
\end{conjecture}

\begin{theorem} Let 
$F_1=T_3^2=
\left[\begin{array}{lll}
0 & 1 & 1\\
1 & 0 & 1\\
1 & 1 & 0
 \end{array}
 \right]$,
and  
$F_2=[T_3^2,T_3^3]=
\left[\begin{array}{llll}
0 & 1 & 1 & 1\\
1 & 0 & 1 & 1\\
1 & 1 & 0 & 1
 \end{array}
 \right]$.
Then $\sat(n,F_1)=\sat(n,F_2)=3n-2$ for any $3\le n\le 6$.
Moreover, for any $n\ge 7$, $\sat(n,F_1)\le 3n-2$ and $\sat(n,F_2)\le 3n-2$ as well.
\end{theorem}

\begin{proof}
Let $M_n=[T_n^0,T_n^1,T_n^n,\tilde{T}_n^2]$, where $\tilde{T}_n^2\subseteq T_n^2$ consists of all those columns of $T_n^2$ which have precisely one entry equal to~1 either in the first or in the $n$th row (but not in both), \di{e.g.}, for $n=5$ we obtain

$$
 M_5=\left[
 \begin{array}{lllllllllllll}
0 & 1 & 0 & 0 & 0 & 0 & 1 & 1 & 1 & 1 & 0 & 0 & 0\\
0 & 0 & 1 & 0 & 0 & 0 & 1 & 1 & 0 & 0 & 1 & 0 & 0\\
0 & 0 & 0 & 1 & 0 & 0 & 1 & 0 & 1 & 0 & 0 & 1 & 0\\
0 & 0 & 0 & 0 & 1 & 0 & 1 & 0 & 0 & 1 & 0 & 0 & 1\\
0 & 0 & 0 & 0 & 0 & 1 & 1 & 0 & 0 & 0 & 1 & 1 & 1
\end{array}
\right].
$$
Clearly, $e(M_n) = e(T_n^0) + e(T_n^1) + e(T_n^n) + e(\tilde{T}_n^2) = 1 + n +
1 + 2n-4 = 3n-2$. Moreover, since $\tilde{T}_n^2$ is $F_1$-admissible we get
that $M_n$ is both $F_1$ and $F_2$ admissible. Now we show that $M_n$ is
$F_1$-saturated. Indeed, pick
any column $C=(c_1,\dots,c_n)^T$ which is not present in $M_n$. Such a column must contain at
least $2$ ones and $1$ zero. Let $1\le i,j,k\le n$ be the indices such
that $c_i=0$, $c_j=c_k=1$. If $i=1$ or $i=n$, then the matrix
$[M_n,C](\{i,j,k\},)$ contains $F_1$. Otherwise, $c_1=c_n=1$, and there also
exists $1<i<n$ such that $c_i=0$. Here $[M_n,C](\{1,i,n\},)$ contains $F_1$. 
Thus $M_n$ is $F_1$ saturated and, since it must contain $T_n^n$ is a column, $M_n$ is also $F_2$-saturated.
We conclude that
$\sat(n,F_1)\le 3n-2$ and $\sat(n,F_2)\le 3n-2$ for any $n\ge 3$. 
A computer search~\cite{dudek+pikhurko+thomason:09:arxiv} yields that these inequalities are
equalities when $n=3,\dots,6$. 
\end{proof}

\begin{conjecture}
Let $F_1=T_3^2$ and $F_2=[T_3^2,T_3^3]$. Then $\sat(n,F_1)=\sat(n,F_2)=3n-2$ for every $n\ge 7$.
\end{conjecture}

\begin{remark}
It is not hard to see that $\sat(n, F_1) \ge n+c\sqrt{n}$ for some absolute
constant~$c$ and all $n\ge 3$. 
Indeed, let $M$ be an $n\times(n+2+\lambda)$ $F_1$-saturated
matrix of size $\sat(n, F_1)$ for some~$\lambda=\lambda(n)$. We may assume
that  $M(,[n+2]) = [T_n^0,T_n^1,T_n^n]$. Suppose that $\lambda\le n$ for
otherwise we are done. 
Moreover, we assume that every column
of matrix $M([\lambda], \{n+3,\dots,n+2+\lambda\})$ contains at least one
entry equal to~1 (trivially, there must be a permutation of the rows of~$M$ satisfying
this requirement). We claim that all rows of $M(\{\lambda+1,\dots,n\},
\{n+3,\dots,n+2+\lambda\})$  are different. Suppose not. Then, there are
indices $\lambda+1\le i,j\le n$ such that $M(i,\{n+3,\dots,n+2+\lambda\}) =
M(j,\{n+3,\dots,n+2+\lambda\})$. Now consider a column $C$ in which the only
nonzero entries correspond to~$i$ and~$j$. Clearly, $C$ is not present in~$M$,
since the first $\lambda$ entries of~$C$ equal~0. Moreover, since~$M$ is
$F_1$-saturated, the matrix $[M, C]$ contains~$F_1$. In other words, there are
three rows in $M$ which form $F_1$ as a submatrix. Note that the $i$th and
$j$th row must be among them. But this is not possible since $F_1$ has no pair
of equal rows.  

Let $M_0 = M(\{\lambda+1,\dots,n\}, \{n+3,\dots,n+2+\lambda\})^T$. Clearly,
$M_0$ is $F_1$-admissible. Anstee and Sali showed (see Theorem~1.3
in~\cite{anstee+sali:05}) that $\forb(\lambda,F_1)=O(\lambda^2)$. That means
that $n-\lambda = O(\lambda^2)$, and consequently,
$\lambda=\Omega(\sqrt{n})$. Hence, $\sat(n,F_1) = e(M) \ge
n+\Omega(\sqrt{n})$, as required.
\end{remark}

\section*{Acknowledgements}

We are grateful to the referees for their careful reading and insightful
comments.

\bibliographystyle{amsplain}

\bibliography{sat-mat-refs}


\clearpage
\section*{C source code of \texttt{satmat}}

{\footnotesize{
\begin{verbatim}
/*  
 *   Copyright (C) 2009 Andrzej Dudek, Oleg Pikhurko and Andrew Thomason
 *
 *   This program is free software: you can redistribute it and/or modify
 *   it under the terms of the GNU General Public License as published by
 *   the Free Software Foundation, either version 3 of the License, or
 *   (at your option) any later version.
 *
 *   This program is distributed in the hope that it will be useful,
 *   but WITHOUT ANY WARRANTY; without even the implied warranty of
 *   MERCHANTABILITY or FITNESS FOR A PARTICULAR PURPOSE.  See the
 *   GNU General Public License for more details.
 *
 *   You should have received a copy of the GNU General Public License
 *   along with this program.  If not, see <http://www.gnu.org/licenses/>.
 *
 *
 *   satmat.c: this program finds the smallest number of columns (saturated 
 *             matrix) with n rows without having an r-row matrix F
 * 
 *   compile: 
 *       % gcc -Wall -ansi -O4 -o satmat satmat.c
 *       
 *   compile w/ debug option: 
 *       % gcc -Wall -ansi -O4 -o satmat -DDEBUG satmat.c
 *       
 *   usage:
 *       % ./satmat    
 * 
 *   --------------------------------------------------------  
 *   Authors: Andrzej Dudek      adudek@andrew.cmu.edu
 *            Oleg Pikhurko      pikhurko@andrew.cmu.edu
 *            Andrew Thomason    A.G.Thomason@dpmms.cam.ac.uk 
 *   --------------------------------------------------------  
 */

#include <stdlib.h>
#include <stdio.h>
#include <time.h>
#include <string.h>
 
#define MAXN 10     /* max value of n */

#define MAXR 5      /* max value of r */

#define MAXNCR 252  /* should be MAXN choose MAXR */

#define MAXS 50     /* max size of strong sat set */

int sub_mask[1 << MAXN][MAXNCR]; /* the sub-column mask indexed by some r rows,
                                    e.g. n=5 r=3  subset 4 in colex is 10011
                                    30th column is 11110  subcolumn = 110 = 6
                                    so subcol[30][4] = ..111011111 */

int subs_left[MAXS][MAXNCR]; /* allowed sub-columns indexed by rows found so
                                far; e.g. if the first 18 columns under
                                consideration contain, in rows indexed by
                                subset 4, sub-columns of types 0,1,4,5
                                then subs_left[18][4] = ~110011 = ~51 */

int poss_cols[(1 << MAXN) * MAXS]; /* list of columns to try next */

int T[MAXR][1 << MAXR];            /* family of T_k^l matrices */

int F[MAXR][1 << MAXR];            /* forbidden matrix */

int colT[MAXR+2];                  /* colT[l] is the first column of T_k^l in T */

int irrel[1 << MAXN];              /* no sub-column is in F */

int mask;                          /* column representation of matrix F */



void gen_T(int r)
{
    int rows[MAXR];
    register int i, l, col;

    /* init T */
    memset((void *)T, 0, MAXR*(1 << MAXR)*sizeof(int));
    
    colT[0] = 0;
    col = 1;
    for (i = 0; i < r; i++) { /* precisely (i+1) 1s per column */
    
        colT[i+1] = col;
    
        for (l = 0; l < i+1; l++) {
            rows[l] = l; 
        }    
        
        do {
 
            for (l = 0; l < i+1; l++) {
                T[rows[l]][col] = 1;
            }    
            col++;

            /* find next subset of i 1s */
            for (l = 0; l < i; l++) {
                if (rows[l] + 1 != rows[l + 1]) {
                    break;
                }    
            }           
                    
            rows[l] += 1;

            while (--l >= 0) {
                rows[l] = l;
            }       
 
        } while (rows[i] < r);
    }
    colT[i+1] = col;
}

void gen_mask(int r, int c)
{
    register int i, j, col;

#ifdef DEBUG
    printf("\n F mask\n");
#endif    
    mask = 0;
    for (j = 0; j < c; j++) {
        for (col = 0, i = 0; i < r; i++) {
            if (F[i][j] == 1) {
                col |= (1<<i);
            }
        }
        mask |= (1 << col);
#ifdef DEBUG        
        printf("  mask[%d] --> %d\n", j, col);  
#endif        
    }

#ifdef DEBUG
    printf(" mask --> %d\n", mask);
#endif    
}

void gen_irrel(int n, int min0, int min1, int *irr)
{
    register int i,j,ones;

#ifdef DEBUG
    printf(" \n min0 --> %d\n", min0);
    printf(" min1 --> %d\n", min1);
#endif

    /* column with less than min1 1s or more than (n-min0) 1s */
    for (*irr = 0, i = 0; i < (1 << n); i++) {
        for (ones = 0, j = 0; j < n; j++) { /* count 1s */
            if (i & (1 << j)) {
            	ones++;
            }
        }
        
        if (ones < min1 || ones > n-min0) {
            irrel[(*irr)++] = i;
        }    
    }    

#ifdef DEBUG
    printf("\n irrelevant columns --> ");    
    for (i = 0; i < *irr; i++) {
        printf("%2d ", irrel[i]); 
    }
    printf("\n");    
#endif    
} 

int is_irrel(int col, int irr)
{
    register int i;

    /* binary search could be more effective */
    for (i = 0; i < irr; i++) {
        if (irrel[i] == col) {
            return 1;
        }
    }

    return 0; 
}

void input(int *n, int *target, int *r, int *c, int *irr)
{
    int i,j,k,col,min0,min1;
    char ans;
    
    printf("\n************************** Wellcome to SATMAT! **************************\n");
    printf("This program checks if sat(n,F) <= target, for a given an r-row matrix F.\n");
    printf("If so, an n-row F-saturated matrix is produced.\n");
    printf("*************************************************************************");
    
    printf("\n input n -->  ");
    scanf("%d", n);
    if (*n < 2  ||  *n > MAXN) {
        printf(" n is not in the range 2 to %d (MAXN)\n\n", MAXN);
        exit(0);
    }

    printf("\n input r -->  ");
    scanf("%d", r);
    if (*r < 2  ||  *r > MAXR) {
        printf(" r is not in the range 2 to %d (MAXR)\n\n", MAXR);
        exit(0);
    }
    if (*r > *n) {
        printf(" r is not in the range 2 to %d \n\n", *n);
        exit(0);
    }

    printf("\n target  -->  ");
    scanf("%d", target);
    if (*target < 2  ||  *target >= MAXS) {
        printf(" target is not in the range 2 to %d (MAXS-1)\n\n", MAXS - 1);
        exit(0);
    }
    
    printf("\n input F as a union of Ti matrices\n\n");
    gen_T(*r);
    for (i = 0; i < *r; i++) {
        for (j = 0; j <= *r; j++) {
            if (i == (*r-1)/2) {
                printf("  T%d=|", j);
            }
            else {
                printf("     |");
            }
        
            for (k = colT[j]; k < colT[j+1]; k++) {
                printf("%d", T[i][k]);
            }
            printf("|");    
        }
        printf("\n"); 
    }
    
    col = 0;
    min0 = min1 = *r;
    for (i = 0; i <= *r; i++) {
        printf("\n add T%d (y/n) --> ", i); 
        scanf(" %c", &ans);
        if (ans == 'y') { /* copy Ti into F */
            for (j = 0; j < *r; j++) {
                for (k = colT[i]; k < colT[i+1]; k++) {
                    F[j][col+k-colT[i]] = T[j][k];
                }
            }
            col += (colT[i+1]-colT[i]);
            
            /* at least this number of 0s in very column*/
            if (*r-i < min0) { 
                min0 = *r-i; 
            }
            
            /* at least this number of 1s in very column*/
            if (i < min1) {
                min1 = i;       
            }
        }
    }
    *c = col;
    
    /* print F */
    printf("\n");
    for (i = 0; i < *r; i++) {
        if (i == (*r-1)/2) {
            printf("  F=|");
        }
        else {
            printf("    |");
        }
        
        for (j = 0; j < col; j++) {
            printf("%d", F[i][j]);
        }
        
        printf("|\n");    
    }
    
    /* find F mask */
    gen_mask(*r, *c);
    
    /* find irrelevan columns */
    gen_irrel(*n, min0, min1, irr);
}

void setup(int n, int r, int c, int *ncr)
{
    int rows[MAXR];
    register int l, col, subr, subcol;

    /* subset of r rows */
    for (l = 0; l < r; l++) {
        rows[l] = l; 
    }
    
    /* index of subset in colex */
    subr = 0;            

    do {
        /*  for each column, extract subcolumn indexed by rows*/
        for (col = 0; col < (1 << n); col++) {
            for (subcol = 0, l = 0; l < r; l++) {
                if (col & (1 << rows[l])) {
                	subcol |= (1 << l);
                }
            }

            sub_mask[col][subr] = (1 << (1 << r)) - 1 - (1 << subcol);
            sub_mask[col][subr] &= mask;
        }
 
        /* find next colex subset of rows */
        for (l = 0; l < r - 1; l++) {
            if (rows[l] + 1 != rows[l + 1]) {
                break;
            }
        }
        rows[l] += 1;
        while (--l >= 0) {
            rows[l] = l;
        }
        subr += 1;
        
    } while (rows[r - 1] < n); /* don't do colex beyond n */

    *ncr = subr;  /* we stopped at subset n choose r */
    fflush(stdout);
}

void print_sol(int *cols, int n, int l, int irr)
{
    int matrix[MAXN][MAXS];
    int i, j;

    /* columns from the solution */
    for (j = 0; j < l; j++) {
        for (i = n-1; i >= 0; i--) {
            matrix[i][j] = ((1<<i) & cols[j]) ? 1 : 0;
        }
    }

    /* irrelevant columns */
    for (j = 0; j < irr; j++) {
        for (i = n-1; i >= 0; i--) {
            matrix[i][l+j] = ((1<<i) & irrel[j]) ? 1 : 0;
        }
    }

    l += j;

    for (i = 0; i < n; i++) {
    
        if (i == (n-1)/2) {
            printf("  M_{%d,%2d}=|", n, l);
        }
        else {         
            printf("           |");
        }
    
        for (j = 0; j < l; j++) {
            printf("%d", matrix[i][j]);
        }
        printf("|\n");
    }
    fflush(stdout);
}

void find_sat(register int n, register int r, register int ncr, register int target, int irr)
{
    int cols[MAXS];                   /* current set of columns */
    int *ptr[MAXS], *to[MAXS];        /* pointers into poss_cols */
    register int l, m, subr, col;
    register int *p, *q, *s, *t, *u;

    register int strong = target + 1; /* smallest strong set so far */
    register int examples = 0;

    int gone_critical[MAXNCR]; /* sub-columns which became critical */
    int *gc[MAXS];             /* pointers into gone_critical */
    int crit_at[MAXNCR];       /* when a sub-column went critical */

    time_t start, stop;
    double diff;
#ifdef DEBUG    
    long long counter = 0;
#endif    

    start = time(NULL);

    for (m = 0, col = 0; col < (1 << n); col++) {
        if (is_irrel(col,irr) == 0) {
            poss_cols[m++] = col;
        }
    }

    ptr[0] = poss_cols;
    to[0] = poss_cols + m;

#ifdef DEBUG
    printf("\n number of relevant columns --> %d\n", m);
    fflush(stdout);
#endif

    /* there is only one relevant column */
    if (m == 1) {
        print_sol((int*)0, n, 0, irr);	
        strong  = 0;
        examples  = 1;
    }

    gc[0] = gone_critical;
    for (m = 0; m < ncr; m++) {
        crit_at[m] = -1;
    }

    l = 0;

    /* we have set of l + 1 columns */
    do { /* main loop */  

#ifdef DEBUG
        counter++;
#endif        

        cols[l] = *ptr[l]++;
 
        if (l + 1 >= strong  ||  ptr[l] > to[l]) { 
            l -= 1; 
            if (l < 0) {
                break;
            }
            goto undo; 
        }    

#ifdef DEBUG
        /* info */
        if (l < 5) {
            for (m = 0; m <= l; m++) {
                printf("%3x", cols[m]);
            }
            putchar('\r');
            fflush(stdout);
        }
#endif

        /* make up subs_left mask after column l */
        if (l == 0) {
            for (subr = 0; subr < ncr; subr++) {
                subs_left[0][subr] = sub_mask[cols[l]][subr];
            }
            p = subs_left[0];
        }  
        else {
            p = subs_left[l - 1]; 
        }

        q = subs_left[l];
        s = sub_mask[cols[l]];
        t = p + ncr;
        u = gc[l];

        while (p < t) {
            m = *q++ = *p++ & *s++;
            if ((m & (m - 1)) == 0) { /* only bad subcolumn */
                m = ncr - 1 - (t - p);
                if (crit_at[m] < 0) { /* subcolumn not critical so far */
                    crit_at[m] = l;
                    *u++ = m;
                }
            }
        }
        gc[l + 1] = u;

        p = ptr[l];  /* find columns which extend this set */
        q = s = to[l];
        t = gc[l + 1];
        while (p < s) {
            m = *p++;
            u = gc[l];  /* check only subcolumns which went critical */
            while (u < t) {
                subr = *u++;
                if ((subs_left[l][subr] & sub_mask[m][subr]) == 0) { 
                    m = 0;  
                    break; 
                }
            }

            if (m) {
                *q++ = m;
            }
        }
    
        ptr[l + 1] = to[l];
        to[l + 1] = q;

        /* if possible, go round again with larger l */
        if (ptr[l + 1] < to[l + 1]) {
            l += 1;  
            continue; 
        }
 
        /* if not, set is inextendible upwards; is it strongly saturated? */
        subr = 0;

        for (col=0; col < cols[0]; col++) {
            if (is_irrel(col,irr) == 1) {
                continue;
            }
            for (subr = 0; subr < ncr; subr++) {
                if ((subs_left[l][subr] & sub_mask[col][subr]) == 0) {
                    break;
                }
            }
            if (subr == ncr) {
                break;
            }
        }
            
        if (subr < ncr) {
            for (m = 0; m < l; m++) {
                for (col = cols[m] + 1; col < cols[m + 1]; col++) {
                    if (is_irrel(col,irr) == 1) {
                        continue;
                    }
                    for (subr = 0; subr < ncr; subr++) {
                        if ((subs_left[l][subr] & sub_mask[col][subr]) == 0) {
                            break;
                        }
                    }

                    if (subr == ncr) {
                        break;
                    }
                }

                if (subr == ncr) {
                    break;
                }
            }
        }

        if (col == cols[l]) { /* we have a strong set size l + 1*/
            if (strong < l + 1) {
                continue; /* not best possible */
            }
            if (strong > l + 1) {
                strong = l + 1; /* best so far */
                examples = 0;
            }
            
            examples += 1;

            if (examples == 1) {

                print_sol(cols, n, strong, irr);

                stop = time(NULL);
                diff = difftime(stop, start);
                printf("\n after --> %.1f min\n", diff/60);
                fflush (stdout); 
            }
        }

        /* now go again with same l */
        undo:
 
        u = gc[l];  /* sub-cols critical at l no longer are */
        t = gc[l + 1];
        while (u < t) {
            crit_at[*u++] = -1;
        }

    } while (l >= 0); /* main loop */

#ifdef DEBUG
    printf("\n number of iterations --> %lld\n", counter);
#endif
    
    stop = time(NULL);
    diff = difftime(stop, start);
    printf("\n time taken --> %.1f min\n", diff/60);
    
    printf("\n******************************* Solution ********************************\n");
    if (examples > 0) {
        printf("%28s sat(%d,F) = %d\n", "",n, strong+irr);
    }
    else {
        printf("%28s sat(%d,F) > %d\n", "", n, target+irr);	
    }
    printf("*************************************************************************\n");
}

int main(int argc, char *argv[])
{
    int n, r, c, ncr, target, irr;

    input(&n, &target, &r, &c, &irr);
    setup(n, r, c, &ncr);
    find_sat(n, r, ncr, target-irr, irr);

    exit(0);
}
\end{verbatim}
}}

\end{document}